\documentclass{article}

\usepackage{amsmath}
\usepackage{amsfonts}
\title{\LARGE \textbf{Sequential edge-coloring on the subset of vertices of almost regular
graphs}}
\author{Petros A. Petrosyan\\ \\Department of Informatics and Applied Mathematics,\\
Yerevan State University, 0025, Armenia\\
Institute for Informatics and Automation Problems,\\
National Academy of Sciences, 0014, Armenia\\ E-mail:
pet\_petros@ipia.sci.am}

\begin{document}
\textheight = 18.1cm \textwidth = 11.4cm \maketitle

\begin{abstract}
Let $G$ be a graph and $R\subseteq V(G)$. A proper edge-coloring of
a graph $G$ with colors $1,\ldots,t$ is called an $R$-sequential
$t$-coloring if the edges incident to each vertex $v\in R$ are
colored by the colors $1,\ldots,d_{G}(v)$, where $d_{G}(v)$ is the
degree of the vertex $v$ in $G$. In this note, we show that if $G$
is a graph with $\Delta(G)-\delta(G)\leq 1$ and
$\chi^{\prime}(G)=\Delta(G)=r$ ($r\geq 3$), then $G$ has an
$R$-sequential $r$-coloring with $\vert R\vert \geq
\left\lceil\frac{(r-1)n_{r}+n}{r}\right\rceil$, where $n=\vert
V(G)\vert$ and $n_{r}=\vert\{v\in V(G):d_{G}(v)=r\}\vert$. As a
corollary, we obtain the following result: if $G$ is a graph with
$\Delta(G)-\delta(G)\leq 1$ and $\chi^{\prime}(G)=\Delta(G)=r$
($r\geq 3$), then $\Sigma^{\prime}(G)\leq \left\lfloor\frac
{2n_{r}(2r-1)+n(r-1)(r^{2}+2r-2)}{4r}\right\rfloor$, where
$\Sigma^{\prime}(G)$ is the edge-chromatic sum of $G$.\\
\end{abstract}

\section{Introduction}

In this note we consider graphs which are finite, undirected, and
have no loops or multiple edges. Let $V(G)$ and $E(G)$ denote the
sets of vertices and edges of a graph $G$, respectively. The degree
of a vertex $v\in V(G)$ is denoted by $d_{G}(v)$ and the chromatic
index of $G$ by $\chi^{\prime }\left(G\right)$. For a graph $G$, let
$\Delta(G)$ and $\delta(G)$ denote the maximum and minimum degrees
of vertices in $G$, respectively. An $(a,b)$-biregular bipartite
graph $G$ is a bipartite graph $G$ with the vertices in one part all
having degree $a$ and the vertices in the other part all having
degree $b$. The terms and concepts that we do not define can be
found in \cite{West}.

A proper edge-coloring of a graph $G$ is a mapping $\alpha:
E(G)\rightarrow \mathbf{N}$ such that $\alpha(e)\neq
\alpha(e^{\prime})$ for every pair of adjacent edges
$e,e^{\prime}\in E(G)$. If $\alpha$ is a proper edge-coloring of a
graph $G$, then $\Sigma^{\prime}(G,\alpha)$ denotes the sum of the
colors of the edges of $G$. For a graph $G$, define the
edge-chromatic sum $\Sigma^{\prime}(G)$ as follows:
$\Sigma^{\prime}(G)=\min_{\alpha}\Sigma^{\prime}(G,\alpha)$, where
minimum is taken among all possible proper edge-colorings of $G$. A
proper $t$-coloring is a proper edge-coloring which makes use of $t$
different colors. If $\alpha $ is a proper $t$-coloring of $G$ and
$v\in V(G)$, then $S\left(v,\alpha\right)$ denotes set of colors
appearing on edges incident to $v$. Let $G$ be a graph and
$R\subseteq V(G)$. A proper edge-coloring of a graph $G$ with colors
$1,\ldots,t$ is called an $R$-sequential $t$-coloring if the edges
incident to each vertex $v\in R$ are colored by the colors
$1,\ldots,d_{G}(v)$.

The concept of sequential edge-coloring of graphs was introduced by
Asratian \cite{Asr}. In \cite{Asr,AsrKam1},
he proved the following result.\\

\noindent\textbf{Theorem 1.} If $G=(X\cup Y,E)$ is a bipartite graph
with $d_{G}(x)\geq d_{G}(y)$ for every $xy\in E(G)$, where $x\in X$
and $y\in Y$, then $G$ has an $X$-sequential $\Delta(G)$-coloring.\\

On the other hand, in \cite{AsrKam1} Asratian and Kamalian showed
that the problem of deciding whether a bipartite graph $G=(X\cup
Y,E)$ with $\Delta(G)=3$ has an $X$-sequential $3$-coloring is
$NP$-complete. Some other results on sequential edge-colorings of
graphs were obtained in \cite{Kam,KamSweden}. In particular, in
\cite{KamSweden} Kamalian proved the following result.\\

\noindent\textbf{Theorem 2.} If $G$ is an $(r-1,r)$-biregular
($r\geq 3$) bipartite graph with $n$ vertices, then $G$ has an
$R$-sequential $r$-coloring with $\vert R\vert \geq
\left\lceil\frac{rn}{2r-1}\right\rceil$.\\

In this note we generalize last theorem. As a corollary, we also
obtain the following result: if $G$ is a graph with
$\Delta(G)-\delta(G)\leq 1$ and $\chi^{\prime}(G)=\Delta(G)=r$
($r\geq 3$), then $\Sigma^{\prime}(G)\leq \left\lfloor\frac
{2n_{r}(2r-1)+n(r-1)(r^{2}+2r-2)}{4r}\right\rfloor$.\\

\section{The Result}

\noindent\textbf{Theorem 3.} If $G$ a graph with
$\Delta(G)-\delta(G)\leq 1$ and $\chi^{\prime}(G)=\Delta(G)=r$
($r\geq 3$), then $G$ has an $R$-sequential $r$-coloring with $\vert
R\vert \geq \left\lceil\frac{(r-1)n_{r}+n}{r}\right\rceil$, where
$n=\vert V(G)\vert$ and $n_{r}=\vert\{v\in V(G):d_{G}(v)=r\}\vert$.

\noindent\textbf{Proof.} Since $\chi^{\prime}(G)=\Delta(G)=r$, there
exists a proper $r$-coloring $\alpha$ of the graph $G$ with colors
$1,2,\ldots,r$. For $i=1,2,\ldots,r$, define the set $V_{\alpha}(i)$
as follows:
\begin{center}
$V_{\alpha}(i)=\left\{v\in V(G):i\notin S(v,\alpha)\right\}$.
\end{center}

Clearly, for any $i^{\prime},i^{\prime\prime}, 1\leq
i^{\prime}<i^{\prime\prime}\leq r$, we have

\begin{center}
$V_{\alpha}(i^{\prime})\cap
V_{\alpha}(i^{\prime\prime})=\emptyset$~~
and~~~$\underset{i=1}{\overset{r}{\bigcup
}}V_{\alpha}(i)=V(G)\setminus V_{r}$,
\end{center}
where $V_{r}=\{v\in V(G):d_{G}(v)=r\}$.\\

Hence,

\begin{center}
$n-n_{r}=\vert V(G)\vert-\vert V_{r}\vert =\left\vert
\underset{i=1}{\overset{r}{\bigcup
}}V_{\alpha}(i)\right\vert=\underset{i=1}{\overset{r}{\sum }}\vert
V_{\alpha}(i)\vert$.
\end{center}

This implies that there exists $i_{0}$, $1\leq i_{0}\leq r$, for
which $\vert V_{\alpha}(i_{0})\vert \geq
\left\lceil\frac{n-n_{r}}{r}\right\rceil$. Let $R=V_{r}\cup
V_{\alpha}(i_{0})$. Clearly, $\vert R\vert\geq n_{r}+
\left\lceil\frac{n-n_{r}}{r}\right\rceil$.

If $i_{0}=r$, then $\alpha$ is an $R$-sequential $r$-coloring of
$G$; otherwise define an edge-coloring $\beta$ as follows: for any
$e\in E(G)$, let

\begin{center}
$\beta(e)=\left\{
\begin{tabular}{ll}
$\alpha(e)$, & if $\alpha(e)\neq i_{0},r$,\\
$i_{0}$, & if $\alpha(e)=r$,\\
$r$, & if $\alpha(e)=i_{0}$.\\
\end{tabular}%
\right.$
\end{center}

It is easy to see that $\beta$ is an $R$-sequential $r$-coloring of
$G$ with $\vert R\vert \geq
\left\lceil\frac{(r-1)n_{r}+n}{r}\right\rceil$. $\square$\\

\noindent\textbf{Corollary 1.} If $G$ is an $(r-1,r)$-biregular
($r\geq 3$) bipartite graph with $n$ vertices, then $G$ has an
$R$-sequential $r$-coloring with $\vert R\vert \geq
\left\lceil\frac{rn}{2r-1}\right\rceil$.\\

\noindent\textbf{Corollary 2.} If $G$ is a graph with
$\Delta(G)-\delta(G)\leq 1$ and $\chi^{\prime}(G)=\Delta(G)=r$
($r\geq 3$), then\
\begin{center}
$\Sigma^{\prime}(G)\leq \left\lfloor\frac
{2n_{r}(2r-1)+n(r-1)(r^{2}+2r-2)}{4r}\right\rfloor$\\
\end{center}
\noindent\textbf{Proof.} Let $\alpha $ be an $R$-sequential
$r$-coloring of $G$ with $\vert R\vert \geq
\left\lceil\frac{(r-1)n_{r}+n}{r}\right\rceil$ described in the
proof of Theorem 3. Now, we have

\begin{eqnarray*}
\Sigma^{\prime}(G)\leq\Sigma^{\prime}\left(G,\alpha\right) &\leq &
\frac{\frac{n_{r}\cdot
r(r+1)}{2}+\left\lceil\frac{n-n_{r}}{r}\right\rceil\frac{r(r-1)}{2}+\left(n-n_{r}-\left\lceil\frac{n-n_{r}}{r}\right\rceil\right)\frac{(r+2)(r-1)}{2}}{2}\\
&\leq &\frac{\frac{n_{r}\cdot
r(r+1)}{2}+\frac{(n-n_{r})r(r-1)}{2r}+\left(n-n_{r}-\frac{n-n_{r}}{r}\right)\frac{(r+2)(r-1)}{2}}{2}\\
&=& \frac{\frac{n_{r}\cdot
r(r+1)}{2}+\frac{(n-n_{r})r(r-1)}{2r}+\frac{(n-n_{r})(r+2)(r-1)^{2}}{2r}}{2}\\
&=& \frac{n_{r}\cdot
r(r+1)}{4}+\frac{(n-n_{r})(r-1)(r^{2}+2r-2)}{4r}\\
&=& \frac{2n_{r}(2r-1)+n(r-1)(r^{2}+2r-2)}{4r}.
\end{eqnarray*}
$\square$\\

\end{document}